\documentclass[11pt]{amsart}
\usepackage{amsmath,amssymb}
\usepackage{amscd,eucal,amsthm}
\newtheorem{Theorem}{Theorem}

\newtheorem{Lemma}{Lemma}

\theoremstyle{remark}

\renewcommand{\to}[1][]{\xrightarrow{#1}}

\renewcommand{\gg}{{\mathfrak{g}}}
\renewcommand{\ll}{{\mathfrak{l}}}

\newcommand{\s}{{\mathfrak{s}}}
\newcommand{\ofr}{{\mathfrak{o}}}
\newcommand{\p}{{\mathfrak{p}}}

\newcommand{\hh}{{\mathfrak{h}}}

\newcommand{\qq}{\mathbb{Q}}
\newcommand{\cc}{\mathbb{C}}

\newcommand{\A}{\mathcal{A}}
 
\DeclareMathOperator{\gr}{gr}

 \DeclareMathOperator{\tr}{tr}

 \DeclareMathOperator{\rk}{rk}

\renewcommand{\phi}{\varphi}

 \makeatletter
\def\@mult#1{\raise #1\rlap{$\cdot$}\lower #1\rlap{$\cdot$}\cdot}
\def\did{\mathrel{\@mult{3pt}}}
\makeatother
\def\openrow#1#2#3{\setbox0=\vbox{\hbox
    {\vrule height#2 width#3\kern#2\vrule height#2 width0pt}\hrule height#3}
    \hbox{\leaders\copy0\hskip#1\wd0\vrule width#3}}
\def\row#1#2#3{\vbox{\hrule height#3\openrow{#1}{#2}{#3}}}
\def\Yr#1{\row{#1}{1.5ex}{.1ex}}

\def\DY#1\endDY{\baselineskip=1ex\lineskip=0pt\lineskiplimit=0pt{\vcenter
    {\Yr#1}}}
\def\openclm#1#2#3{\setbox0=\vbox{\hrule height#3\hbox
    {\vrule width0pt\kern#2\vrule width#3 height#2}}\vtop
    {\leaders\copy0\vskip#1\ht0\hrule height#3}}
\def\clm#1#2#3{\hbox{\vrule width#3\openclm{#1}{#2}{#3}}}
\def\Yc#1{\clm{#1}{1.5ex}{.1ex}}

\def\CDY#1\endCDY{{\vcenter{\hbox{\Yc#1}}}}

\author{L.~G.~Rybnikov}
\title{Centralizers of certain quadratic elements in Poisson--Lie algebras and Argument Shift method}
\address{Poncelet laboratory (Independent University of Moscow and CNRS) and Moscow State University,
department of Mechanics and Mathematics}
\email{leo.rybnikov@gmail.com}
\thanks{The work was partially supported by RFBR grant 04-01-00702, RFBR grant 05
01 00988-a and RFBR grant 05-01-02805-CNRSL-a.}

\begin{document}
\maketitle

\section{Introduction.} Let $\gg$ be a semisimple complex Lie algebra.
The universal enveloping algebra $U(\gg)$ bears a natural
filtration by the degree with respect to the generators. The
associated graded algebra $\gr U(\gg)$ is naturally isomorphic to
the symmetric algebra $S(\gg)=\cc[\gg^*]$ by the
Poincar\'e--Birkhoff--Witt theorem. The commutator operation on
$U(\gg)$ defines a Poisson bracket on $S(\gg)$, which we call the
\emph{Poisson--Lie bracket}.

The \emph{argument shift method} gives a way to construct
Poisson-commutative subalgebras in $S(\gg)$. The method is as
follows. Let $ZS(\gg)=S(\gg)^{\gg}$ be the center of $S(\gg)$ with
respect to the Poisson bracket, and let $\mu\in\gg^*$ be a regular
semisimple element. Then the algebra $A_{\mu}\subset S(\gg)$
generated by the elements $\partial_{\mu}^n\Phi$, where $\Phi\in
ZS(\gg)$, (or, equivalently, generated by central elements of
$S(\gg)=\cc[\gg^*]$ shifted by $t\mu$ for all $t\in\cc$) is
Poisson-commutative and has maximal possible transcendence degree
equal to $\frac{1}{2}(\dim\gg+\rk\gg)$ (see \cite{MF}). Moreover,
the subalgebras $A_{\mu}$ are maximal Poisson-commutative
subalgebras in $S(\gg)$ (see \cite{T}). In \cite{V}, the
subalgebras $A_{\mu}\subset S(\gg)$ are named the {\em
Mischenko--Fomenko subalgebras}.

Let $\hh\subset\gg$ be a Cartan subalgebra of the Lie algebra
$\gg$. We denote by $\Delta$ and $\Delta_+$ the root system of
$\gg$ and the set of positive roots, respectively. Let
$\alpha_1,\dots,\alpha_l$ be the simple roots. Fix a
non-degenerate invariant scalar product $(\cdot,\cdot)$ on $\gg$
and choose from each root space $\gg_{\alpha},\ \alpha\in\Delta,$
a nonzero element $e_{\alpha}$ such that
$(e_{\alpha},e_{-\alpha})=1$. Set
$h_{\alpha}:=[e_{\alpha},e_{-\alpha}]$, then for any $h\in\hh$ we
have $(h_{\alpha},h)=\langle\alpha,h\rangle$. The elements
$e_{\alpha}\ (\alpha\in\Delta)$ together with
$h_{\alpha_1},\dots,h_{\alpha_l}\in\hh$ form a basis of $\gg$.

We identify $\gg$ with $\gg^*$ via the scalar product
$(\cdot,\cdot)$ and assume that $\mu$ is a regular semisimple
element of the fixed Cartan subalgebra $\hh\subset\gg=\gg^*$. The
linear and quadratic part of the Mischenko--Fomenko subalgebras
can be described as follows \cite{F}:
\begin{gather*}A_\mu\cap\gg=\hh,\\ A_\mu\cap S^2(\gg)=S^2(\hh)\oplus
Q_\mu,\ \text{where}\
Q_\mu=\{\sum\limits_{\alpha\in\Delta_+}\frac{\langle\alpha,h\rangle}{\langle\alpha,\mu\rangle}e_{\alpha}e_{-\alpha}|h\in\hh\}.\end{gather*}

The main result of the present paper is the following

\begin{Theorem}\label{main} For generic $\mu\in\hh$
(i.e. for $\mu$ in the complement to a certain countable union of
Zariski-closed subsets in  $\hh$), the algebra $A_\mu$ is the
Poisson centralizer of the subspace $Q_\mu$ in $S(\gg)$.
\end{Theorem}

In \cite{ChT,NO,T2} the Mischenko--Fomenko subalgebras were lifted
(quantized) to the universal enveloping algebra, i.e. the family
of commutative subalgebras $\A_\mu\subset U(\gg)$ such that
$\gr\A_\mu=A_\mu$ was constructed for any classical Lie algebra
$\gg$ (i.e. $\s\ll_r$, $\s\ofr_r$, $\s\p_{2r}$). In \cite{R} we do
this (by different methods) for any semisimple $\gg$.

We deduce the following assertion from Theorem~\ref{main}.

\begin{Theorem}\label{lift} For generic $\mu\in\hh$
there exist no more than one commutative subalgebra $\A_\mu\subset
U(\gg)$ satisfying $\gr\A_\mu=A_\mu$.
\end{Theorem}

This means that there is a \emph{unique} quantization of
Mischenko--Fomenko subalgebras. In particular, the methods of
\cite{ChT,NO,T2} and \cite{R} give the same for classical Lie
algebras. In the case $\gg=\gg\ll_n$ the assertion of
Theorem~\ref{lift} was proved by A.~Tarasov \cite{T3} for any
regular $\mu\in\hh$.

I thank E.~B.~Vinberg for attention to my work and useful
discussions.

\section{Proof of Theorem~\ref{main}} Note that the set $E_n\subset\hh$ of such $\mu\in\hh$
that the Poisson centralizer of the space $Q_\mu$ in $S^n(\gg)$
has the dimension greater than  $\dim A_\mu\cap S^n(\gg)$ is
Zariski-closed in $\hh$ for any $n$. Therefore it suffices to
prove that $E_n\ne\hh$ for any $n$. Thus, it suffices to prove the
existence of $\mu\in\hh$ satisfying the conditions of the Theorem.
\begin{Lemma}\label{irrational}
There exist $\mu,h\in\hh$ such that numbers
$\frac{\langle\alpha,h\rangle}{\langle\alpha,\mu\rangle}\
(\alpha\in\Delta_+)$ are linearly independent over $\qq$.
\end{Lemma}
\begin{proof}
Choose $\mu$ such that the values $\alpha_i(\mu)$ are
algebraically independent over $\qq$ for simple roots $\alpha_i$.
Since there are no proportional positive roots, the numbers
$\frac{1}{\langle\alpha,\mu\rangle},\ \alpha\in\Delta_+,$ are
linearly independent over $\qq$. Choose $h$ such that the values
$\langle\alpha,h\rangle$ are nonzero rational numbers. Then the
numbers $\frac{\langle\alpha,h\rangle}{\langle\alpha,\mu\rangle},\
\alpha\in\Delta_+,$ are linearly independent over $\qq$.
\end{proof}

Choose $\gamma\in\gg^*$ such that $\gamma(h_{\alpha_i})=1$ for any
simple root $\alpha_i$ and $\gamma(e_{\alpha})=0$ for
$\alpha\in\Delta$. We define a new Poisson bracket
$\{\cdot,\cdot\}_{\gamma}$ on $S(\gg)$ by setting
$\{x,y\}_{\gamma}=\gamma([x,y])$ for $x,y\in\gg$. This bracket is
compatible with the Poisson--Lie bracket, i.e. the linear
combination $t\{\cdot,\cdot\}+(1-t)\{\cdot,\cdot\}_{\gamma}$ is a
Poisson bracket on $S(\gg)$ (i.e. satisfies the Jacobi identity)
for any $t\in\cc$. Moreover, for $t\ne0$, the corresponding
Poisson algebras are isomorphic. Namely, denote by $S(\gg)_t$ the
algebra $S(\gg)$ equipped with the Poisson bracket
$t\{\cdot,\cdot\}+(1-t)\{\cdot,\cdot\}_{\gamma}$; then for $t\ne0$
the Poisson algebra isomorphism $\psi_t:S(\gg)_1\to S(\gg)_t$ is
defined on the generators $x\in\gg$ as follows:
$\psi_t(x)=t^{-1}x+t^{-2}(1-t)\gamma(x)$. Clearly, we have
$\psi_t(Q_\mu)=Q_\mu$.

\begin{Lemma}\label{lem_tr_deg}
The transcendence degree of the Poisson centralizer of the
subspace $Q_\mu$ in $S(\gg)_0$ is not greater than
$\frac{1}{2}(\dim\gg+\rk\gg)$ for some $\mu\in\hh$.
\end{Lemma}
\begin{proof}
Choose $\mu$ and $h$ as in Lemma~\ref{irrational} and set
$q=\sum\limits_{\alpha\in\Delta_+}\frac{\langle\alpha,h\rangle}{\langle\alpha,\mu\rangle}e_{\alpha}e_{-\alpha}\in
Q_\mu$. For any $f\in S(\gg)$, we have
$\{q,f\}_{\gamma}=\sum\limits_{\alpha\in\Delta_+}\gamma(h_{\alpha})\frac{\langle\alpha,h\rangle}{\langle\alpha,\mu\rangle}(e_{-\alpha}\frac{\partial
f}{\partial e_{-\alpha}}-e_{\alpha}\frac{\partial f}{\partial
e_{\alpha}})$. In particular,
\begin{multline}\label{monom}
\{q,\prod\limits_{i=1}^{l}h_{\alpha_i}^{m_i}\prod\limits_{\alpha\in\Delta_+}
e_{\alpha}^{n_{\alpha}}e_{-\alpha}^{n_{-\alpha}}\}_{\gamma}=\\=
\sum\limits_{\alpha\in\Delta_+}\gamma(h_{\alpha})\frac{\langle\alpha,h\rangle}{\langle\alpha,\mu\rangle}(n_{-\alpha}-n_{\alpha})\prod\limits_{i=1}^{l}h_{\alpha_i}^{m_i}\prod\limits_{\alpha\in\Delta_+}
e_{\alpha}^{n_{\alpha}}e_{-\alpha}^{n_{-\alpha}}.
\end{multline}
For any $\alpha=\sum\limits_{i=1}^{l}k_i\alpha_i\in\Delta_+$, we
have
$\gamma(h_{\alpha})=\sum\limits_{i=1}^{l}k_i\in\qq\backslash\{0\}$.
Since the numbers
$\frac{\langle\alpha,h\rangle}{\langle\alpha,\mu\rangle}$ are
linearly independent over $\qq$, the right hand part
of~(\ref{monom}) is zero iff $n_{\alpha}-n_{-\alpha}=0$ for any
$\alpha\in\Delta_+$. This means that the Poisson centralizer of
$q$ in $S(\gg)_0$ is linearly generated by monomials having equal
degrees in $e_{\alpha}$ and $e_{-\alpha}$ for any
$\alpha\in\Delta_+$, i.e. the Poisson centralizer of $q$ in
$S(\gg)_0$ is generated (as a commutative algebra) by the elements
$h_{\alpha_i}\ (i=1,\dots, l)$ and $e_{\alpha}e_{-\alpha}\
(\alpha\in\Delta_+)$. Therefore, the transcendence degree of the
Poisson centralizer of $q$ in $S(\gg)_0$ is equal to
$\frac{1}{2}(\dim\gg+\rk\gg)$.
\end{proof}

By Lemma~\ref{lem_tr_deg}, the transcendence degree of the Poisson
centralizer of the subspace $Q_\mu$ in $S(\gg)_t$ is not greater
than $\frac{1}{2}(\dim\gg+\rk\gg)$ for generic $t$. Since the
Poisson algebras $S(\gg)_t$ are isomorphic to each other for
$t\ne0$, this lower bound of the transcendence degree holds for
any $t\in\cc$. Let $Z\subset S(\gg)$ be the Poisson centralizer of
$Q_\mu$ in $S(\gg)_1$. Since $\tr\ \deg(Z)\le\tr\ \deg(A_\mu)$ and
$A_\mu\subset Z$, we see that each element of $Z$ is algebraic
over $A_\mu$. By Tarasov's results \cite{T}, the subalgebra
$A_\mu$ is algebraically closed in $S(\gg)_1$, hence, $Z=A_\mu$.
Theorem~\ref{main} is proved.

\section{Proof of Theorem~\ref{lift}} By \cite{V}, the subspace $A_\mu^{(2)}=\cc+\hh+S^2(\hh)+Q_\mu\subset S(\gg)^{(2)}$
can be uniquely lifted to a commutative subspace
$\A_\mu^{(2)}\subset U(\gg)^{(2)}$ (this subspace is the image of
$A_\mu^{(2)}$ under the symmetrization map). By
Theorem~\ref{main}, any lifting $\A_\mu\subset U(\gg)$ of $A_\mu$
is the centralizer of the subspace $\A_\mu^{(2)}$ in $U(\gg)$ in
$U(\gg)$ for generic $\mu$. Theorem~\ref{lift} is proved.

\end{document}